\documentclass[12pt]{article}

\usepackage{amsfonts}
\usepackage{graphicx}
\usepackage{enumerate}
\usepackage[shortlabels]{enumitem}
\usepackage{mathtools}
\usepackage{hyperref}
\usepackage{color}
\usepackage{xcolor}
\usepackage{fullpage}
\usepackage{tikz}
\usepackage{cite}  

\usepackage{amsmath, amssymb,latexsym}
\usepackage{soul}
\usepackage{float}
\usepackage{comment}

\def\F{\mathbb{F}}

\def\K{\mathbb{K}}

\def\N{\mathbb{N}}
\def\Q{\mathbb{Q}}
\def\Z{\mathbb{Z}}

\def\OK{\mathcal{O}_\mathbb{K}}

\DeclareMathOperator{\lcm}{lcm}

\newtheorem{theorem}{Theorem}[section]
\newtheorem{proposition}[theorem]{Proposition}

\newtheorem{lemma}[theorem]{Lemma}

\newtheorem{remark}[theorem]{{\sc Remark}}

\usepackage[normalem]{ulem}
\definecolor{brilliantrose}{rgb}{1.0, 0.33, 0.64}
\definecolor{myviolet}{rgb}{0.21, 0.0, 0.85}
\definecolor{amethyst}{rgb}{0.6, 0.4, 0.8}
\definecolor{carrotorange}{rgb}{0.93, 0.57, 0.13}
\definecolor{cutepink}{rgb}{1.0, 0.2, 0.6}

\title{Eventual periodicity of the Smith forms of integer matrix powers}
\author{Vanni Noferini\thanks{Aalto University, Department of Mathematics and Systems Analysis, P.O. Box 11100, FI-00076, Aalto, Finland. Supported by a Research Council of Finland grant (decision number 370932). Email: vanni.noferini@aalto.fi}}
\date{11 June 2026}

\begin{document}

\maketitle

\begin{abstract}
    We prove that the Smith forms of the powers of an integer square matrix behave in an eventually periodic manner. More precisely, if $\mathrm{SF}(M)$ denotes the Smith form of $M \in \Z^{m \times m}$, then for every $A \in \Z^{m \times m}$ there exist $n_0 \in \N$, an integer $T \geq 1$, and a constant diagonal matrix $D \in \Z^{m \times m}$ such that $n \geq n_0$ implies $\mathrm{SF}(A^{n+T})=D \cdot \mathrm{SF}(A^n)$. This provides an eventually affirmative answer to a conjecture posed in 2013 by R. Bruner.  We also show that both $n_0$ and $T$ can be arbitrarily large. 
\end{abstract}


\bigskip


\section{Introduction}

The goal of this paper is to prove the following theorem.

\begin{theorem}\label{thm1}
    Let $A \in \Z^{m \times m}$ and define the sequence $(S_n)_{n \in \N}$, where $S_n$ is the Smith form of $A^n$ normalized to have nonnegative diagonal entries. Then, there exists an eventually periodic sequence $(D_n)_{n \in \N} \subset \Z^{m \times m}$ of diagonal matrices such that $S_{n+1}=D_n S_n$ for all $n \in \N$.
\end{theorem}

Theorem \ref{thm1} was conjectured by R. Bruner \cite{mo}, in the stronger form where $(D_n)_n$ is periodic (when restricted to $n \geq 1$) and not just eventually periodic. However, the stronger form does not hold, and in fact the threshold $n_0$ before $D_n$ starts behaving periodically can be arbitrarily large. We provide a family of counterexamples in  Proposition \ref{prop:no} (whose proof is postponed to Subsection \ref{sec:smallproof}), showing also that $n_0$ can grow linearly with $m$.

\begin{proposition}\label{prop:no}
    If, in Theorem \ref{thm1}, we replace ``eventually periodic" with ``periodic", then the statement is false.
\end{proposition}

In spite of the promising comments to \cite{mo}, to our knowledge a proof to Theorem \ref{thm1} was thus far lacking. For a start, it is interesting to observe that the period in Theorem \ref{thm1} is unbounded even for bounded size $m$. For example, let $p \in \N$ be a prime and observe that the following $2 \times 2$ example has period $p$:
\[ A = \begin{bmatrix}
    p & 1\\
    0 & p
\end{bmatrix} \Rightarrow A^n = \begin{bmatrix}
    p^n & n p^{n-1}\\
    0 & p^n
\end{bmatrix} \Rightarrow S_n = \begin{cases}
   \begin{bmatrix}
       p^n & 0\\
       0 & p^n
   \end{bmatrix} \ & \ \mathrm{if} \ p \mid n ;\\
   \\
  \begin{bmatrix}
      p^{n-1} & 0\\
      0 & p^{n+1}
  \end{bmatrix} \ & \ \mathrm{otherwise}.
\end{cases} \]
It is clear that similar examples can be constructed for every size, e.g., $A = \begin{bmatrix}
    p&1\\
    0&p
\end{bmatrix} \oplus I_{m-2}$.

There are several reasons why Theorem \ref{thm1} is interesting. From a matrix theoretical viewpoint, there are known and achievable bounds \cite{qs} for the invariant factors of $AB$ in terms of those of $A$ and $B$, where $A$ and $B$ are square matrices of the same size. These bounds are tightest when $A$ and $B$ have coprime determinants, as in that case the Smith form of the product is the product of the Smith forms \cite[Theorem II.15]{newman}. However, {as observed in \cite{mo}}, in some sense the case $B=A^{n-1}$ lies at the opposite corner with respect to \cite[Theorem II.15]{newman}, since $\det(A^{n-1})=\det(A)^{n-1}$ is very far from being coprime with $\det(A)$, unless $A$ is a unimodular matrix. (Of course, a unimodular $A$ is a completely trivial case for Theorem \ref{thm1}, since it implies $(S_n)_n \equiv I$). Hence, the fact that the Smith forms of $A^n$ behave very regularly is surprising. Theorem \ref{thm1} is also motivated by algebraic topology \cite{bruneremail,brunerbook}. For example, if $A \in \Z^{m \times m}$ is the companion matrix of $1-(1-x)^m$, which arises from K-theory \cite{brunerbook}, then the Smith form of $A^n$ describes the structure of the $(2n-1)$-th connective $K$-homology of the cyclic group of order $m$ \cite{bruneremail}.

The manuscript is structured as follows. In Section \ref{sec:bg}, we review preliminary material about periodic sequences, matrix theory, and $p$-adic analysis. We then prove Theorem \ref{thm1} in Section \ref{sec:proof}. Our proof relies on some basic facts in algebraic number theory, but is otherwise elementary. All the number-theoretical results that we need can be found in  \cite[Chapter II]{neukirch}. To make the paper self-contained, they will also be recalled in Subsection \ref{sec:padic}.

\section{Background}\label{sec:bg}

\subsection{Periodic sequences}

A sequence $f : \N \rightarrow \Q$ is called periodic if there exists $0< T \in \N$ such that, {for all $n \in \N$}, $f(n)=f(n+T)$, and eventually periodic if there exist $n_0 \in \N$ and $0 < T \in \N$ such that $n \geq n_0 \Rightarrow f(n+T)=f(n)$. 

\begin{proposition}\label{prop:easy}
    Let $\{f_i(n)\}_{i=1}^k$ be finitely many periodic (resp. eventually periodic) sequences, and define their sum, product and minimum as $\sigma(n)=\sum_{i=1}^k f_i(n)$, $\pi(n)=\prod_{i=1}^k f_i(n)$, $\mu(n)=\min_{1 \leq i \leq k} f_i(n)$. Then, $\sigma(n),\pi(n),\mu(n)$ are all periodic (resp. eventually periodic) sequences. If $f_2(n)$ is nonzero (resp. eventually nonzero), the same holds for the {pointwise} quotient $f_1(n)/f_2(n)$ {(possibly defined only for sufficiently large  $n$)}.
\end{proposition}
\begin{proof}
    Straightforward, taking period $T=\lcm_i T_i$ in all cases.
\end{proof}

\begin{proposition}\label{prop:lesseasy}
Let $f_1(n),\dots,f_k(n)$ be finitely many sequences whose first differences, $\Delta f_i(n):=f_i(n+1)-f_i(n)$, are eventually periodic. Define $g(n)=\min_i f_i(n)$; then, $\Delta g(n):=g(n+1)-g(n)$ is eventually periodic.
\end{proposition}
\begin{proof}
Denote, respectively, the eventual period and threshold of $\Delta f_i$ by $T_i$ and $(n_0)_i$. Define $T:=\lcm T_i$ and $n_0 = \max_i (n_0)_i$. Then, for all $i$ and all $n \geq n_0$, we have $\Delta f_i(n+T)=\Delta f_i(n)$. Hence, defining $s_i:=\sum_{j=n}^{n+T-1} \Delta f_i(j)$, and noting that $s_i$ depends on $i$ but not on $n$,  we have $f_i(n+T)=f_i(n)+ s_i$. Let $S=\{ i : s_i$ is minimal $\}$. Then, for all sufficiently large $n \geq N \geq n_0$, it holds  $g(n)=f_i(n)$ for some $i \in S$. Thus, there exist a pair (depending on $n$) of indices $i,j \in S$ such that
\[ \Delta g(n+T) = f_i(n+T+1)-f_j(n+T) = f_i(n+1) -f_j(n) = \Delta g(n).\]
This shows that $\Delta g$ is eventually periodic with period $T$.
\end{proof}

\subsection{Matrix theory}

We recall that a square matrix over a principal ideal domain (PID) $R$ is called \emph{unimodular} if its determinant is a unit of $R$. Theorem \ref{thm:smith} was proved by H. J. S. Smith \cite{smith} for $R=\Z$ and later by I. Kaplansky \cite{kaplansky} for any elementary divisor domain (a class of rings that includes in particular every PID).

\begin{theorem}[Smith's Theorem, \cite{kaplansky,smith}]\label{thm:smith}
  Let $R$ be a PID and $M \in R^{m \times n}$. There {exist} unimodular matrices $U \in R^{m \times m}, V \in R^{n \times n}$ such that $M=USV$, where $S$ is diagonal and such that $S_{ii}$ divides $S_{i+1,i+1}$, for all values of $i$ for which this makes sense. The matrix $S$, which is uniquely determined by $M$ up to multiplying its diagonal elements by units of $R$, is called the \emph{Smith form} of $M$; its diagonal elements are the \emph{invariant factors} of $M$. In addition, again up to units of $R$ and denoting by $\gamma_i(M)$ the $i$-th \emph{determinantal divisor} of $M$, i.e., $\gamma_i$ is the $\gcd$ of all $i \times i$ minors of $M$ (and $\gamma_0(M)=1$), {whenever $\gamma_{i-1} \neq 0$} we have the formula
  \begin{equation}\label{eq:detdiv}
      S_{ii} = \frac{\gamma_i}{\gamma_{i-1}}.   
  \end{equation}  
\end{theorem}

Let $A \in \Z^{m \times m}$. Recall from \cite[Section 0.8.1]{hj} that the $r$-th \emph{compound matrix} of $A$, denoted by $C_r(A) \in \Z^{\binom{m}{r} \times \binom{m}{r}}$, is the matrix whose entries are all the $r \times r$ minors of $A$. By the Cauchy-Binet formula \cite[Section 0.8.7]{hj}, it is easy to verify that $C_r(A^n)=C_r(A)^n$, for all $r,n$; see also \cite[Equation (0.8.1.1)]{hj}. 

Observe now that Theorem \ref{thm1} is obviously true when $A$ is nilpotent. Indeed, in that case $S_n$ is eventually zero, and therefore we can take $(D_n)_n$ eventually constant. Therefore, we henceforth assume that $A$ is not nilpotent. For every integer matrix $M$, denote by $\gcd(M)$ the greatest common divisor of the entries of $M$, normalized to be nonnegative. In view of \eqref{eq:detdiv} and the properties of the compound matrices outlined above, Theorem \ref{thm2} below effectively describes the behavior of all the invariant factors of the powers of any integral matrix. Indeed, let $\Gamma_n$ be the diagonal matrix whose diagonal entries are the determinantal divisors of $A^n$; by Theorem \ref{thm2} applied to compound matrices {(except the nilpotent ones, which are a trivial special case)}, $\Gamma_{n+1}=E_n \Gamma_n$, where $E_n$ is an eventually periodic sequence of diagonal matrices. This implies by Proposition \ref{prop:easy} that $S_{n+1}=D_n S_n$, and $D_n$ is diagonal and eventually periodic (its diagonal elements are ratios of consecutive diagonal entries of $E_n$, whenever well defined, and can be set to $0$ otherwise). Moreover, recalling that (for all $i$) the $i$-th invariant factor of $A^n$ divides the $i$-th invariant factor of $A^{n+1}$ \cite{qs}, the matrices $D_n$ in Theorem \ref{thm1} must be (or can be taken to be, when $S_n$ has zero diagonal entries) integral. Hence, Theorem \ref{thm1} is implied by Theorem \ref{thm2}.

\begin{theorem}\label{thm2}
    Let $A \in \Z^{m \times m}$ be not nilpotent. The sequence \[ g(n) := \frac{\gcd(A^{n+1})}{\gcd(A^n)}\] is eventually periodic.
\end{theorem}

Given a prime $p \in \N$, recall that the $p$-adic valuation of $n \in \Z$ is defined as
\[ \nu_p(n) = \max \{ e \in \N \ s.t. \ p^e \mid n \}. \]
(Note that $\nu_p(0)=+\infty$.) We now extend the $p$-adic valuation to integer matrices by defining $\nu_p(M):=\nu_p(\gcd(M))$; more details on this idea will be given in Subsection \ref{sec:padic}.

Denote by $S$ the Smith form of $M \in \Z^{m \times m}$ over $\Z$ and (for every prime $p \in \N$) by $S_{(p)}$ the Smith form of $M$ over the localization of $\Z$ at the prime ideal $\langle p \rangle$, $\Z_{(p)} = \{ \frac{a}{b} \ \mid \ \gcd(a,b)=1, \ p \nmid b \}$. Since $\mathrm{frac}(\Z)=\mathrm{frac}(\Z_{(p)}) = \Q$, the rank $r$ of $M$ does not depend on the base ring; in particular, {$S_{(p)}$} and $S$ have the same number of zero diagonal elements. Hence, by \cite[Theorem II.15]{newman} we have $S=\prod_{p} S_{(p)}$, and (recalling that we are not interested in the case $M=S=0$) the product can be taken over the finite set of the primes that divide $\gamma_r(M)$, i.e., the largest nonzero determinantal divisor of $M$. It is thus clear that Theorem \ref{thm2} is in turn implied by Theorem \ref{thm3}, by taking as period $T$ for Theorem \ref{thm2} the least common multiple (over all primes $p$ that divide {the product of the nonzero eigenvalues of $A$, counted with multiplicity; this product is up to sign equal to a coefficient of the characteristic polynomial of $A$, and hence an integer}) of the period $T_p$ for Theorem \ref{thm3}. 

\begin{theorem}\label{thm3}
    Let $A \in \Z^{m \times m}$ be not nilpotent. For every prime $p \in \N$, there exist a constant $0 \leq a \in \Q$ and an eventually periodic sequence $h(n)$ such that
    \[ \nu_p(A^n) = a n + h(n).  \]
    In particular,  $\nu_p(A^{n+1}) - \nu_p(A^n) = a + h(n+1)-h(n)$    is eventually periodic.
\end{theorem}
We defer the proof of Theorem \ref{thm3}, and therefore of Theorem \ref{thm1}, to Section \ref{sec:proof}.

\subsection{Proof of Proposition \ref{prop:no}}\label{sec:smallproof}

Let us take a small break and show that Theorem \ref{thm1} cannot be strengthened to take the original form of Bruner's conjecture \cite{mo}. 

    \begin{proof}[Proof of Proposition \ref{prop:no}]
    {Let $m \geq 3$ and c}onsider the counterexample $A=C \oplus 2$ where $C \in \Z^{(m-1) \times (m-1)}$ is the companion matrix of the polynomial $x^{m-1}-4^m$. For example, when $m=4$ then
    \[  A = \begin{bmatrix}
        0 & 0 & 256 & 0\\
        1 & 0 & 0 & 0\\
        0 & 1 & 0 & 0\\
        0 & 0 & 0 & 2
    \end{bmatrix} \in \Z^{4 \times 4}.  \]
    Let $n=(m-1)k+r$ be the Euclidean division of $n$ by $m-1$, with $k \in \N$ and $0 \leq r \leq m-2$. Then, it is easy to see that the invariant factors of $A^n$ are, up to reordering, $2^n$ (repeated one  time), $4^{mk+m}$ (repeated $r$ times), $4^{mk}$ (repeated $m-1-r$ times). Since $A$ is invertible over $\Q$, there is a unique choice of $(D_n)_n$ that satisfies $S_{n+1}=D_n S_n$ for all $n$, and $(D_n)_n$ cannot be periodic.
    
    Indeed, if $k=0$ (or, equivalently, $n < m-1$), the smallest invariant factors of $C^n$ are $1$. Thus, $\gcd(A^n)= 1$. On the other hand, if $k \geq 1$, every entry of $C^n$ is divisible by $4^{mk}$. Since $2mk - n = (m+1)k - r > 0$, this implies $\gcd(A^n) = 2^n$. This switch implies that the entry $(D_n)_{11}$ changes behavior, causing a break in periodicity.
\end{proof}

\begin{remark}
    In the proof above, the sequence $(D_n)_n$ is nevertheless eventually periodic, and thus compatible with Theorem \ref{thm1}. More precisely, we have $D_{n+T}=D_n$ for all $n \geq n_0$, with $T=n_0=m-1$.
\end{remark}

\subsection{$p$-adic analysis}\label{sec:padic}

Fix a prime $p \in \N$. The $p$-adic valuation extends to $\Q$ as a function $\nu_p: \Q \rightarrow \Z \cup \{ +\infty \}$ in the following way: If $0\neq \frac{a}{b}=q \in \Q$, then $\nu_p(\frac{a}{b})=\nu_p(a)-\nu_p(b).$ The topological completion of $\Q$ with respect to the $p$-adic absolute value, $|q|_p := p^{-\nu_p(q)}$, is the field of \emph{$p$-adic numbers}, denoted by $\Q_p$. Since $\Z \subsetneq \Q \subsetneq \Q_p$, every square integer matrix $A$ is also a matrix with entries in $\Q_p$. However, $\Q_p$ is not algebraically closed, and therefore the $p$-adic eigenvalues of $A$ generally lie in some finite field extension $\K \supseteq \Q_p$. Every such field extension has the following useful properties \cite[Chapter II]{neukirch}:

\begin{enumerate}
    \item \cite[p. 118-120 and p. 150]{neukirch} The $p$-adic valuation extends to $\K$, as a function $\nu_p: \K \rightarrow \Q \cup \{+ \infty \}$ that coincides with $\nu_p$ as defined above for elements of $\Q$, and retains the following properties:
    \begin{itemize}
        \item[(i)] For $x \in \K$, $\nu_p(x) = +\infty \Leftrightarrow x=0$;
        \item[(ii)] For all $a,b \in \K$, $\nu_p(ab)=\nu_p(a)+\nu_p(b)$;
        \item[(iii)] For all $a,b \in \K$, $\nu_p(a+b) \geq \min \{ \nu_p(a),\nu_p(b)  \}$, with equality if $\nu_p(a) \neq \nu_p(b)$.
    \end{itemize}
    \item \cite[Chapter II, Proposition 3.8]{neukirch} The subset
    \[ \mathcal{O}_\K = \{ x \in \K \ s.t. \ \nu_p(x)\geq 0 \} \]
    is a ring called the \emph{ring of integers of $\K$}, having field of fractions $\K$, group of units
    \[ \mathcal{O}^\times_\K = \{ x \in \K \ s.t. \ \nu_p(x) = 0 \},  \]
    and maximal ideal
    \[ \mathfrak{m} = \{ x \in \K \ s.t. \ \nu_p(x) > 0 \}.\]
    \item \cite[p. 150]{neukirch} There exists a positive integer $f$, called \emph{inertia degree}, such that the residue field $\mathcal{O}_\K/ \mathfrak{m}$ is isomorphic to the finite field with $p^f$ elements:
    \[  \mathcal{O}_\K / \mathfrak{m} \cong \F_{p^f}.  \]
    \item \cite[p. 121 and p. 150]{neukirch} There exists an element $\pi \in \mathcal{O}_\K$, called \emph{uniformizer}, such that $\mathfrak{m}=\langle \pi \rangle$ and a positive integer $e$, called \emph{ramification index} and such that $\langle p \rangle = \langle \pi^e \rangle$.
    \item \cite[Chapter II, Proposition 3.9]{neukirch} Generally, every ideal of $\mathcal{O}_\K$ has the form $\langle \pi^L \rangle$, $L \in \N$, and the ring $\OK / \langle \pi^L \rangle$ is finite and has cardinality $p^{Lf}$.
    \item By the above properties, it follows that $\nu_p(\cdot)$ is a discrete valuation on $\K$ and, on $\K^\times=\K \setminus \{0\}$, it takes values in $\frac1e \Z \subsetneq \Q$.
\end{enumerate}

\begin{theorem}[Fermat's Little Theorem]\label{thm:fermat}
    Let $\theta \in \OK^\times$, i.e., $\nu_p(\theta)=0$, let $f$ be the inertia degree, and let $\pi \in \OK$ be a uniformizer. Then, defining $t:=p^f-1$,
    \[  \theta^{t} \equiv 1 \bmod \langle \pi \rangle.   \]
\end{theorem}
\begin{proof}
It follows by \cite[Chapter II, Proposition 3.10]{neukirch}.
\end{proof}

The $p$-adic valuation can be extended to a matrix $M \in \K^{m \times n}$, by taking the entrywise minimum of the scalar $p$-adic valuations \cite{rp2}:
\[ \nu_p(M) := \min_{\scriptsize\begin{matrix}
    1 \leq i \leq m\\
    1 \leq j \leq n
\end{matrix}} \nu_p(M_{ij}).  \]
Note that, if $M$ has entries in $\Z$, this definition is equivalent to $\nu_p(M):=\nu_p(\gcd_\Z(M))$. The matrixwise $p$-adic valuation retains properties (i) and (iii) as in the scalar case \cite[Proposition 3.3]{rp2}, but property (ii) becomes weaker. Indeed, when the product $AB$ is defined, then generally $\nu_p(AB) \geq \nu_p(A) + \nu_p(B)$ \cite[Proposition 3.3]{rp2}. However, equality holds when either $A$ (or $B$, or both) is unimodular, i.e., a unit of $\OK^{m \times m}$ for some $m$ \cite[Proposition 3.3]{rp2}; observe that being a unit of $\OK^{m \times m}$ is equivalent to $\nu_p(A)=\nu_p(\det A)=0$.

\begin{lemma}\label{lem:jordan}
    Let $A \in \Z^{m \times m}$. There exists an invertible matrix $S \in { \mathrm{GL}(m,\K)}$ such that, for all $n \in \N$, $A^nS=SJ^n$ and $J \in \OK^{m \times m}$ is in Jordan canonical form.
\end{lemma}
\begin{proof}
Observe that $\OK$ is integrally closed \cite[p. 121]{neukirch}; since $\Z \subsetneq \OK$, the eigenvalues of $A$ (being the roots of a monic polynomial in $\Z[x]$) must lie in $\OK$, i.e., they have nonnegative $p$-adic valuations. Hence, all the eigenvectors and Jordan vectors of $A$ can be taken to have entries in $\K=\mathrm{frac}(\OK)$. Defining $S$ as the matrix whose columns are such vectors, properly ordered, we conclude by standard linear algebra over the field $\K$ that $S$ is invertible and $A=SJS^{-1}$. The statement then follows easily by induction on $n$.
\end{proof}

\section{Proof of Theorem \ref{thm3}}\label{sec:proof}

\begin{lemma}\label{lem:polyval}
Let ${ 0 \leq } c  \in \Q$ and let $q(x) \in \K[x]$. The sequence $f(n)=\min \{ 0, \nu_p(q(n)) - c \}$ is periodic.
\end{lemma}
\begin{proof}
Let $\theta \in \OK$ be  such that $Q(x):=\theta q(x) \in \OK[x]$; the existence of $\theta$ is guaranteed because $\K$ is the field of fractions of $\OK$.  Moreover, let $d:= c + \nu_p(\theta) { \geq 0} $ and $D:=\lceil d \rceil$. Suppose first that $\nu_p(q(n)) \geq c$, implying  $f(n)=0$. Then, $\nu_p(Q(n)) \geq d$ and, by a formal Taylor expansion, $Q(n+p^D) \equiv Q(n) \bmod \langle p^D \rangle $. Thus, \[ \nu_p (Q(n+p^D)) \geq \min \{ \nu_p(Q(n)), \nu_p(Q(n+p^D)-Q(n))  \} \geq \min \{d,D \}=d,\] and hence
\[  \nu_p(q(n+p^D)) = -\nu_p(\theta) + \nu_p( Q(n+p^D)) \geq  c \Rightarrow f(n+p^D)=0=f(n). \]    
Now suppose instead $\nu_p(q(n)) < c$, implying $\nu_p(Q(n)) < d$. Then,
\[ \nu_p(q(n+p^D)) = -\nu_p(\theta) + \nu_p(Q(n+p^D)) = -\nu_p(\theta) + \nu_p(Q(n)) < c, \]
and hence  $f(n+p^D)=f(n)=\nu_p(q(n)) -c < 0$.
\end{proof}

Inspecting carefully the structure of the proof of Lemma \ref{lem:polyval}, we see that we can make similar statements for a larger class of functions than just polynomials. We record this fact in Theorem \ref{thm:funval}.

\begin{theorem}\label{thm:funval}
   Let ${Q} : \N \rightarrow \OK$ be a sequence and let $n_0,D \in \N$, and $1 \leq T \in \N$ be such that, for all $n \geq n_0$,
   \[ Q(n+T) \equiv Q(n) \bmod \langle p^D \rangle. \]
   Then, for every $c \in \Q$ and every $\theta \in \OK$ satisfying $c + \nu_p(\theta) \leq D$, the sequence $f(n)=\min \{0,\nu_p(\theta^{-1} Q(n))-c \}$ is eventually periodic, with period dividing $T$.
\end{theorem}
\begin{proof}
  Let $n \geq n_0$ and $\Delta_T(n):=Q(n+T)-Q(n)$. Then, $\nu_p(\Delta_T(n)) \geq D$, and hence
  \begin{equation}\label{eq:magic}
      \nu_p(Q(n+T)) \geq \min \{\nu_p(Q(n)), \nu_p(\Delta_T(n))\} \geq \min \{ \nu_p(Q(n)), D\}.
  \end{equation}
If $\nu_p(Q(n)) \geq c + \nu_p(\theta)$ then \eqref{eq:magic} yields $\nu_p(Q(n+T)) \geq c+ \nu_p(\theta)$, and hence $f(n)=f(n+T)=0$.     On the other hand, if $\nu_p(Q(n))<c+\nu_p(\theta) \leq D$, then we can replace the first $\geq$ in \eqref{eq:magic} by $=$, and
further conclude that $\nu_p(Q(n+T)) = \nu_p(Q(n))$ implying $f(n)=f(n+T) < 0.$  
\end{proof}

\begin{theorem}\label{thm:jordanisperiodic}
    Let $J \in \OK^{m \times m}$ be in Jordan canonical form {and not nilpotent}. Then, the sequence $f(n)=\nu_p(J^{n+1})-\nu_p(J^n)$ is eventually periodic.
\end{theorem}

\begin{proof}
{Consider first a single Jordan block $B$ in $J$, with eigenvalue $\lambda \in \OK$. The statement is obvious if the size of the block is $1$. If $\lambda = 0$, then $\nu_p(B^n)=+\infty$ for sufficiently large $n$. Since $J$ is not globally nilpotent, its nilpotent Jordan blocks do not affect the overall minimum $\nu_p(J^n)$ and can be ignored. So let us suppose that the size is $s+1 \geq 2$ and the eigenvalue is $\lambda \neq 0$. Then, $B^n$ is an upper triangular Toeplitz matrix,}
\[ {B^n=}  \begin{bmatrix}
\lambda & 1 &  &\\
& \ddots & \ddots & \\
&& \lambda & 1\\
&&&\lambda
\end{bmatrix}^n = \begin{bmatrix}
\lambda^n & \binom{n}{1} \lambda^{n-1} & \dots & \binom{n}{s} \lambda^{n-s}\\
&\ddots&\ddots&&\\
\end{bmatrix} {\in \OK^{(s+1) \times (s+1)}}. \]
Hence $${\nu_p(B^n)} = n \nu_p(\lambda) + \min_{0 \leq k \leq s} \left\{\nu_p \left( \binom{n}{k} \right) - k \nu_p(\lambda)\right\} = n \nu_p(\lambda)  + \min_{1 \leq k \leq s} \left\{  \min \{ 0, \nu_p \left(  \binom{n}{k}\right) - k \nu_p(\lambda)   \}  \right\}   .$$
Inspecting this formula, we see  that the inner minima are periodic by Lemma \ref{lem:polyval}, since $\binom{n}{k}$ is a polynomial in $n$ with coefficients in $\Q \subseteq \K$ whereas $k \nu_p(\lambda) {\geq 0}$ is a rational constant, independent of $n$. Therefore, {$\nu_p(B^n)-n \nu_p(\lambda)$} is the minimum of {$s$} periodic functions, and hence periodic by Proposition \ref{prop:easy}. It follows that
{$\nu_p(B^{n+1})-\nu_p(B^n)$} is the sum of the constant $\nu_p(\lambda)$ and two periodic functions, and hence periodic by Proposition \ref{prop:easy}.

More generally, if {$J = \bigoplus_{i=1}^d B_i$} is the direct sum of $d \geq 1$ Jordan blocks, then {$\nu_p(J^n)=\min_{1 \leq i \leq d} \nu_p(B_i^n)$}. Hence, we apply Proposition \ref{prop:lesseasy} to conclude that $\nu_p(J^{n+1})-\nu_p(J^n)$ is eventually periodic. 
\end{proof}

\begin{remark}\label{rem:ark}
    If we order the eigenvalues $\lambda_i$ in $J$ so that 
    \[  \nu_p(\lambda_1) \leq \nu_p(\lambda_2) \leq \dots \leq \nu_p(\lambda_d),\]
    combining the proofs of Proposition \ref{prop:lesseasy} and of Theorem \ref{thm:jordanisperiodic} we obtain the formula
    \[ \nu_p(J^n) = n \nu_p(\lambda_1) + f(n) \]
    where $f(n)$ is an eventually periodic sequence.
\end{remark}

\begin{lemma}\label{lem:daje}
   Let $L$ be a positive integer. For all $k=1,\dots,p^L$, it holds
    \[  \nu_p\left( \binom{p^L}{k} \right) = L - \nu_p(k). \]
    Furthermore, if $e$ is a positive integer, there exists $M \in \N$ such that $\nu_p(k)-\frac{k}{e} \leq M$ for all $1 \leq k \in \N$.
\end{lemma}
\begin{proof}
    Clearly,
    \[ \frac{k}{p^L} \binom{p^L}{k} =  \binom{p^L-1}{k-1} . \]
    By Kummer's theorem \cite{kummer}, the $p$-adic valuation of the right hand side is the number of carries in the addition $(k-1)+(p^L-k)=p^L-1$, performed by writing the digits of each integer in base $p$. But such expansion for $p^L-1$ consists of $L$ digits all equal to $p-1$, and hence there are no carries. The first part of the statement follows.

    For the second part, observe that $\nu_p(k) - \frac{k}{e} \leq \log_p(k) - \frac{k}{e} \rightarrow - \infty$ when $k \rightarrow + \infty$. Hence, such a sequence must have a maximum $\hat{M}$, and we can take $M=\lceil \max\{0,\hat{M}\} \rceil$.
\end{proof}

We finally have all the ingredients ready to cook the main dish.

\begin{proof}[Proof of Theorem \ref{thm3}.] By Lemma \ref{lem:jordan}, we have $A^n = S J^n S^{-1}$. Hence, by \cite[Proposition 3.3]{rp2}, there is a constant $0 \leq k:=-\nu_p(S) -\nu_p(S^{-1})  \in \Q$ such that
\[ \nu_p(J^n)-k \leq \nu_p(A^n) \leq \nu_p(J^n)+k. \]
Furthermore, $(A^n)_{ij}$ is a $\K$-linear combination of every element of $J^n$, and hence it has the form (using the notation of Remark \ref{rem:ark}; note that $\lambda_1  \neq 0$ because $A$ is not nilpotent by assumption)
\[   (A^n)_{ij} = \sum_{\ell=1}^d \lambda_\ell^n q_{ij\ell}(n) = \lambda_1^n \left( \sum_{\ell=1}^d q_{ij\ell}(n)\left(\frac{\lambda_\ell}{\lambda_1}\right)^n  \right),  \]
where $q_{ij\ell}(x) \in \K[x]$ and $\theta_\ell:=\lambda_\ell \lambda_1^{-1} \in \OK$. Write now
\[ \alpha_{ij}^{-1}E_{ij}(n) :=\sum_{\ell=1}^d q_{ij\ell}(n) \theta_\ell^n  . \]
{Here, $\alpha_{ij} \in \OK$ is such that $Q_{ij\ell}(n):=\alpha_{ij}q_{ij\ell}(n) \in \OK[n]$ for all $1 \leq \ell \leq d$.  Note that this implies that $E_{ij}(n)$ is a function from $\N$ to $\OK$, being a finite sum of polynomials with coefficients in $\OK$ and exponentials with base in $\OK$.}

We claim that there exists $1 \leq T,D,n_0 \in \N$ such that, for all $i,j$, $E_{ij}(n)$ satisfies the assumptions of Theorem \ref{thm:funval}. Let $F:={\max_n} f(n)$, where $f(n)$ is the eventually periodic sequence defined in Remark \ref{rem:ark}, and note that $F$ is finite because every eventually periodic sequence is bounded. Defining $c:=F+k$, we then have 
\[ c \geq   \nu_p(A^n)- n \nu_p(\lambda_1) = \min_{i,j} \nu_p(\alpha^{-1}_{ij} E_{ij}{(n)}) = \]
\[=\min_{i,j} \min \{ c,\nu_p(\alpha_{ij}^{-1}E_{ij}(n)) \} = c  + \min_{i,j}  \min \{ 0, \nu_p( \alpha_{ij}^{-1}E_{ij}(n)) -c\}    ,\]
{where the inner $\min$ in the last two expressions denotes the ordinary minimum between two numbers taken pointwise for each $(i,j,n)$, and we used the fact (proved in the previous steps) that the outer minimum over $(i,j)$ is bounded above by $c$, so capping
each term at $c$ does not alter the overall minimum. Hence}, $\nu_p(A^n)=n \nu_p(\lambda_1)+h(n)$ where, by Theorem \ref{thm:funval} and Proposition \ref{prop:easy}, $h(n)$ is eventually periodic provided that $D \geq c + \nu_p(\alpha_{ij})$ for all $i,j$. This proves Theorem \ref{thm3}.

It remains to prove the claim, with an appropriate value of $D$. Let $s=\max \{ i : \nu_p(\theta_i)=0\}$, and write
\[ E_{ij}(n) = \sum_{\ell=1}^s Q_{ij\ell}(n) \theta_\ell^n + \sum_{\ell=s+1}^d Q_{ij\ell}(n) \theta_\ell^n. \]
Fix now a nonnegative integer $D \geq c + \max_{i,j}\nu_p(\alpha_{ij})$, and denote by $f$ the inertia degree and by $e$ the ramification index (see Subsection \ref{sec:padic}). Take $n \geq n_0 \geq D \cdot (\min_{\ell > s} \nu_p(\theta_\ell) )^{-1}$ (or $n_0=0$ if $s=d$), let $L= D+M$, where $M$ is the constant defined in Lemma \ref{lem:daje} taking $e$ to be the ramification index. Finally, define $t:=p^f-1 \geq 1$ and $T:=t p^L$. Let us analyze the difference $E_{ij}(n+T)-E_{ij}(n)$. This is the sum over $\ell$ of the addends (all lying in $\OK$)
\begin{equation}\label{eq:addends}
    \theta_\ell^n \left(  \theta_\ell^{T} Q_{ij\ell}(n+T) - Q_{ij\ell}(n) \right).   
\end{equation} 
  If $\ell > s$, then $\nu_p(\theta_\ell^n) \geq n_0 \nu_p(\theta_\ell) \geq D$. Since \eqref{eq:addends} is a multiple of $\theta_\ell^n$, it is congruent to $0$ modulo $\langle p^D \rangle$. When $\ell \leq s$, by Fermat's Little Theorem \ref{thm:fermat}, we have $\theta_\ell^t \equiv 1 \bmod \langle \pi \rangle$, where $\pi$ is a uniformizer. Hence, for some $\zeta \in \OK$, 
\[  \theta^T_\ell -1 = (1 + \zeta \pi)^{p^L} -1 = \sum_{k=1}^{p^L} \binom{p^L}{k}(\zeta \pi)^k. \]
By Lemma \ref{lem:daje}, this implies
\[  \nu_p(\theta_\ell^T-1) \geq \min_{1 \leq k \leq p^L} \left\{ L - \nu_p(k) + \frac{k}{e} \right\} \geq L - M {=} D.   \]
Hence, $\theta_\ell^T \equiv 1 \bmod \langle p^D \rangle$. On the other hand, $Q_{ij\ell}(n+T) - Q_{ij\ell}(n) \equiv 0 \bmod \langle p^D \rangle$ by a formal Taylor expansion. Combining these two congruences, the term between brackets in \eqref{eq:addends} vanishes modulo $\langle p^D \rangle$. This proves the claim.
    
\end{proof}

\section{Acknowledgements}

This paper answers a question that was posed to me by Robert Bruner \cite{bruneremail}, to whom I am obviously indebted. In part, I also benefited from observations made by MathOverflow users in \cite{mo}. I would also like to thank David Carlip, who performed some preliminary investigations on the problem during a summer research internship, and Giulio Peruginelli, with whom I discussed the problem. {Finally, I am grateful to an anonymous reviewer who shared some remarks that improved the presentation.}

\end{document}